\documentclass[11pt,letterpaper]{amsart}

\newcommand{\1}{{1\!\!1}}
\usepackage{amssymb}
\usepackage{xspace}
\usepackage{pb-diagram,lamsarrow,pb-lams}

\usepackage{array}
\newcolumntype{C}{>{$}c<{$}}
\newcolumntype{L}{>{$}l<{$}}
\newcolumntype{R}{>{$}r<{$}}

\makeatletter

\let\c@equation=\c@subsection

\makeatother

\def\to{\mathchoice
{\longrightarrow}
{\rightarrow}
{\rightarrow}
{\rightarrow}}

\def\mapsto{\mathchoice
{\DOTSB\mapstochar\longrightarrow}
{\DOTSB\mapstochar\rightarrow}
{\DOTSB\mapstochar\rightarrow}
{\DOTSB\mapstochar\rightarrow}}

\def\To{\mathchoice
{\Longrightarrow}
{\Rightarrow}
{\Rightarrow}
{\Rightarrow}}

\def\hookrightarrow{\mathchoice
{\DOTSB\lhook\joinrel\relbar\joinrel\rightarrow}
{\DOTSB\lhook\joinrel\rightarrow}
{\DOTSB\lhook\joinrel\rightarrow}
{\DOTSB\lhook\joinrel\rightarrow}}

\newtheorem{theorem}[subsection]{Theorem}
\newtheorem*{theorem*}{Theorem}
\newtheorem{proposition}[subsection]{Proposition}
\newtheorem{lemma}[subsection]{Lemma}
\newtheorem{corollary}[subsection]{Corollary}

\theoremstyle{definition}

\newtheorem{definition}[subsection]{Definition}
\newtheorem{example}[subsection]{Example}
\newtheorem{remark}[subsection]{Remark}

\newcommand{\Z}{\mathbb{Z}}
\newcommand{\C}{\mathbb{C}}

\newcommand{\Q}{\mathbb{Q}}

\newcommand{\om}{\omega}
\newcommand{\eps}{\varepsilon}
\renewcommand{\o}{\otimes}
\newcommand{\ohat}{\Hat{\otimes}}

\DeclareMathOperator{\Hom}{Hom}
\DeclareMathOperator{\Aut}{Aut}

\DeclareMathOperator{\U}{U}

\DeclareMathOperator{\Ind}{Ind}
\DeclareMathOperator{\Res}{res}

\newcommand{\p}{\partial}
\newcommand{\Om}{\Omega}
\newcommand{\CM}{\mathcal{M}}
\newcommand{\CE}{\mathcal{E}}
\DeclareMathOperator{\gr}{gr}
\renewcommand{\*}{\cdot}

\newcommand{\Mbar}{\overline{\mathcal{M}}}
\renewcommand{\]}{{]\!]}}
\renewcommand{\[}{{[\![}}
\DeclareMathOperator{\ch}{ch}
\newcommand{\Dual}{\vee}

\DeclareMathOperator{\Tr}{Tr}
\DeclareMathOperator{\Exp}{Exp}
\DeclareMathOperator{\Log}{Log}
\newcommand{\bull}{\bullet}
\renewcommand{\SS}{\mathbb{S}}
\newcommand{\point}{\Spec(\C)}

\DeclareMathOperator{\K}{\mathsf{K}}
\DeclareMathOperator{\D}{\mathsf{D}}
\DeclareMathOperator{\DD}{\mathsf{T}}
\newcommand{\MHM}{{\mathsf{MHM}}}
\newcommand{\RR}{{\mathsf{R}}}
\newcommand{\TT}{{\mathbb{T}}}
\DeclareMathOperator{\Serre}{\mathsf{e}}
\newcommand{\Cat}{{\mathsf{Cat}}}

\newcommand{\Ab}{{\mathsf{Ab}}}

\newcommand{\Proj}{{\mathsf{Proj}}}
\newcommand{\Var}{{\mathsf{Var}}}

\newcommand{\FF}{\mathsf{F}}
\newcommand{\G}{{\mathbb{G}}}
\newcommand{\LL}{\mathsf{L}}
\newcommand{\SSS}{\mathsf{S}}
\newcommand{\op}{\circ}

\newcommand{\CR}{\mathcal{R}}
\newcommand{\CC}{\mathcal{C}}

\newcommand{\n}{\mathbf{n}}
\DeclareMathOperator{\SL}{SL}
\DeclareMathOperator{\SU}{SU}
\DeclareMathOperator{\GL}{GL}
\DeclareMathOperator{\Spec}{Spec}

\def\({(\!(}
\def\){)\!)}
\newcommand{\Id}{\text{id}}

\renewcommand{\a}{\mathbf{a}}
\renewcommand{\k}{\mathsf{k}}
\renewcommand{\t}{\mathbf{t}}

\renewcommand{\S}{\mathsf{S}}
\newcommand{\s}{\mathsf{s}}
\renewcommand{\l}{\mathsf{L}}
\newcommand{\HH}{\mathsf{H}}
\newcommand{\CF}{\mathcal{F}}

\newcommand{\CU}{\mathcal{U}}

\newcommand{\jbar}{\bar\jmath}
\renewcommand{\a}{\mathbf{a}}

\newcommand{\Cech}{\v{C}ech\xspace}
\newcommand{\Resolve}{\mathcal{L}}

\DeclareMathOperator{\real}{\mathsf{real}}
\DeclareMathOperator{\Ob}{Ob}
\DeclareMathOperator{\Cone}{\mathsf{cone}}
\DeclareMathOperator{\Tot}{Tot}

\newcommand{\embedding}{immersion\xspace}
\newcommand{\embeddings}{immersions\xspace}

\newcommand{\polynomial}{characteristic\xspace}
\newcommand{\polynomials}{characteristics\xspace}

\begin{document}

\title{Resolving mixed Hodge modules on configuration spaces}

\author{E. Getzler}

\address{Max-Planck-Institut f\"ur Mathematik, Gottfried-Claren-Str.\ 26,
D-53225 Bonn, Germany}

\curraddr{Department of Mathematics, Northwestern University, Evanston, IL
60208-2730, USA}

\email{getzler@math.nwu.edu}

\maketitle

If $\pi:X\to S$ is a continuous map of locally compact topological spaces
and $n$ is a natural number, let $X^n/S$ be the $n$th fibred power of $X$
with itself, and let $\FF(X/S,n)$ be the configuration space, whose fibre
$\FF(X/S,n)_s$ over a point $s\in S$ is a configuration of $n$ distinct
points in the fibre $X_s$, or equivalently, the complement of the
$\binom{n}{2}$ diagonals in $X^n/S$. Let $j(n):\FF(X/S,n)\hookrightarrow
X^n/S$ be the natural open embedding.

An essential r\^ole is played in this paper by the higher direct image with
compact support $f_!$ (this is often written $\mathsf{R}f_!$) associated to
a continuous map $f:X\to Y$. This is a functor from the derived category of
sheaves of abelian groups on $X$ to the derived category of sheaves of
abelian groups on $Y$. For example, if $Y$ is a point and $\CF$ is a sheaf
on $X$, then $f_!\CF$ is a complex of abelian groups whose cohomology is
$H^\bull_c(X,\CF)$. If $f$ is either a closed or an open embedding, $f_!$
takes sheaves on $X$ to sheaves on $Y$. (See Section 1 of Verdier
\cite{Verdier2}.)

Given a sheaf $\CF$ of abelian groups on $X^n/S$, we introduce in this
paper a natural resolution $\Resolve^\bull(X/S,\CF,n)$ of the sheaf
$j(n)_!j(n)^*\CF$ by sums of terms of the form $i(J)_!i(J)^*\CF$. (Here,
$i(J)$ is the closed embedding of a diagonal in $X^n/S$.) This resolution
has the property that if $\CF$ is an $\SS_n$-equivariant sheaf (where the
symmetric group $\SS_n$ acts on $X^n/S$ by permuting the factors in the
fibred product), the resolution is $\SS_n$-equivariant as well. For
example, if $n=2$, we have the exact sequence of sheaves
\begin{equation} \label{diagonal}
0 \to j(2)_!j(2)^*\CF \to \CF \to i_!i^*\CF \to 0 ,
\end{equation}
where $i:X\hookrightarrow X^2/S$ is the diagonal embedding. When
$X$ is a Riemann surface, this resolution was introduced by Looijenga
\cite{Looijenga}.

Let $\pi(n):\FF(X/S,n)\to S$ and $\pi(n):X^n/S\to S$ be the projections to
$S$. (We denote them by the same symbol, since confusion is hardly likely
to arise.) The objects $\pi(n)_!j(n)^*\CF$ and
$\pi(n)_!\Resolve^\bull(X/S,\CF,n)$ are isomorphic in the derived category
of sheaves on $S$. We use this isomorphism to calculate the
$\SS_n$-equivariant Euler characteristic of $\pi(n)_!j(n)^*\CF$.

Our main result is that a similar resolution exists when $\pi:X\to S$ is a
quasi-projective morphism of complex varieties and $\CF$ is a mixed Hodge
module: in this way, we obtain a new proof of the formula of \cite{I} for
the Serre \polynomial of the configuration space $\FF(X,n)$, with the
virtue that it applies with no modification to the relative case. (The
Serre \polynomial is the Euler characteristic of $H^\bull_c(\FF(X,n),\Q)$
in the Grothendieck group of mixed Hodge structures.)

A similar spectral sequence has been obtained by Totaro
\cite{Totaro}. Since he works with cohomology, and not cohomology with
compact support, his results depend on the relative dimension of $X/S$ and
require the projection $\pi:X\to S$ to be smooth; however, in that case,
our spectral sequence is equivalent to his.

To extend our resolution from sheaves to mixed Hodge modules, we have to
modify it, since $i^*$ is not a $t$-exact functor of mixed Hodge modules
(or of perverse sheaves, which underly mixed Hodge modules), even when $i$
is a closed \embedding; already for $n=2$, it is well-known that
\eqref{diagonal} must be replaced by an exact triangle. This difficulty is
overcome by introducing \Cech resolutions for the sheaves
$i(J)_!i(J)^*\CF$, which are constructed using the property of mixed Hodge
modules that $f_!$ is $t$-exact for open affine \embeddings. We actually
show that the proof works whenever we have three of Grothendieck's six
operations, namely $\o$, $f_!$ and $f^*$, satisfying an analogue of the
above condition on open affine \embeddings.

If we apply our resolution to the universal elliptic curve, we obtain a
formula for the relative $\SS_n$-equivariant Serre \polynomial of
$\CM_{1,n}/\CM_{1,1}$. Eichler-Shimura theory, which calculates the
cohomology of polynomials of the Hodge local system on $\CM_{1,1}$, then
leads to a formula for the $\SS_n$-equivariant Serre \polynomial of
$\CM_{1,n}$.

The first value of $n$ for which $\CM_{1,n}$ has a contribution from the
cusp forms, and hence has non-Tate cohomology, is $n=11$, for which the
(non-equivariant) Serre \polynomial may be calculated by \eqref{non-equi} to
be
\begin{multline*}
\Serre(\CM_{1,11}) =
\LL^{11}-330\,\LL^9+4575\,\LL^8-30657\,\LL^7+124992\,\LL^6 - \SSS_{12} \\
-336820\,\LL^5+584550\,\LL^8-406769\,\LL^3-865316\,\LL^2+2437776\,\LL
-1814400 ;
\end{multline*}
here, $\LL$ denotes the mixed Hodge structure $\Q(-1)$, and $\SSS_{12}$ is
a two-dimensional Hodge structure of weight $11$, associated to the
discriminant cusp form $\Delta$. (We do not reproduce the equivariant Serre
\polynomial $\Serre^{\SS_{11}}(\CM_{1,11})$ for lack of space: there are
$56$ irreducible representations of $\SS_{11}$, although not all of these
occur.) In particular, $\CM_{1,11}$ has Euler characteristic $-302400$.

The virtual Euler characteristic of the orbifold $\CM_{1,1}$ equals
$-1/12$. (This is a special case of the formula of Harer and Zagier
\cite{HZ}, but is easy to prove directly, using the standard fundamental
domain for the action of $\SL(2,\Z)$ on the upper half-plane.) It follows
by induction on $n$, using the fibrations $\CM_{1,n}\to\CM_{1,n-1}$, that
the virtual Euler characteristic of $\CM_{1,n}$ equals
$(-1)^n(n-1)!/12$. For $n\ge5$, $\CM_{1,n}$ is a fine moduli space (that
is, no automorphism of an elliptic curve fixes $5$ points), and thus its
virtual Euler characteristic equals its Euler characteristic. The agreement
between the resulting formula for $\chi(\CM_{1,11})$ and the value which we
have calculated provides a (modest) consistency check between our work and
that of Harer and Zagier. We show in Proposition \ref{don} that our formula
for the Serre \polynomial of $\CM_{1,n}$ does give the correct value of
$\chi(\CM_{1,n})$, for all $n\ge5$.

In a sequel to this paper \cite{semi}, we show how to sum the Serre
\polynomials of the strata of $\Mbar_{1,n}$ to obtain a formula for its
Hodge polynomial. For example,
\begin{multline*}
\Serre(\Mbar_{1,11}) = \LL^{11} + 2037\,\LL^{10} + 213677\,\LL^9 +
4577630\,\LL^8 + 30215924\,\LL^7 + 74269967\,\LL^6 \\ {} - \SSS_{12} +
74269967\,\LL^5 + 30215924\,\LL^4 + 4577630\,\LL^3 + 213677\,\LL^2 +
2037\,\LL + 1 .
\end{multline*}

\subsection*{Outline of the paper} In Section 1, we explain the relationship
between Arnold's calculation of the cohomology of the configuration spaces
$\FF(\C,n)$ and the theory of Stirling numbers of the first and second
kinds.

Section 2 is devoted to the construction of the resolution in the simpler
case of sheaves of abelian groups. This is generalized in Section 3 to the
cases of perverse sheaves and of mixed Hodge modules.

In Section 4, we apply the associated spectral sequence to generalize the
formula of \cite{I} for the $\SS_n$-equivariant Serre \polynomial of
$\FF(X,n)$ to the relative case.

In Section 5, we apply the formulas of Section 4 to calculate the
$\SS_n$-equivariant Serre \polynomial of the moduli space $\CM_{1,n}$.

\subsection*{Acknowledgments}

I wish to thank the Department of Mathematics at the Universit\'e de
Paris-VII, MIT and the Max-Planck-Institut f\"ur Mathematik in Bonn for
their hospitality during the inception, elaboration and completion of this
paper, respectively. I am grateful to E. Looijenga for introducing me to
the Eichler-Shimura theory, to D. Zagier for his help with the proof of
Proposition \ref{don}, and to the anonymous referee for a number of
excellent suggestions and corrections.

The author is partially supported by a research grant of the NSF and a
fellowship of the A.P. Sloan Foundation.

\section{The combinatorics of partitions and Stirling numbers}

In this section, we recall the cohomology of the configuration space
$\FF(\C,n)$ and its relationship with the Stirling numbers. We give more
detail on the theory of Stirling numbers than is necessary, since it
illuminates the combinatorics which we will apply to construct our
resolution.

\subsection{Partitions}
A partition $J$ of $n$ is a decomposition of the set $\{1,\dots,n\}$ into
disjoint non-empty subsets: for example, the partitions of $\{1,2,3,4\}$
are
\begin{gather*}
\{1,2,3,4\}\quad\{12,3,4\}\quad\{13,2,4\}\quad\{14,2,3\}\quad\{23,1,4\}
\quad\{24,1,3\}\quad\{34,1,2\} \\
\{123,4\}\quad\{124,3\}\quad\{134,2\}\quad\{234,1\}\quad\{12,34\}\quad\{13,24\}
\quad\{14,23\}\quad\{1234\} ,
\end{gather*}
where we abbreviate the subset $\{i_1,\dots,i_\ell\}$ to $i_1\dots
i_\ell$. We denote the subsets of $J$ by $\{J_1,\dots,J_k\}$, in no
particular order. Denote by $\S(n,k)$ the set of partitions of $n$ into $k$
non-empty subsets.

Associated to a partition $J$ of $n$ is an equivalence relation on
$\{1,\dots,n\}$, such that $i\sim_Jj$ iff $i$ and $j$ lie in the same part
of $J$. The set of all partitions of $n$ is a poset: if $J$ and $K$ are
partitions, $J\prec K$ iff $i\sim_Jj$ implies that $i\sim_Kj$, that is, iff
$K$ is coarser than $J$.

If $\a=(a_n\mid n\ge1)$ is a sequence of natural numbers, let
$|\a|=\sum_{n=1}^\infty na_n$.
\begin{lemma} \label{partitions}
The exponential generating function of the number $p(\a)$ of partitions of
$|\a|$ into $a_j$ subsets of size $j$, $j\ge1$, is
$$
B(\t) = \sum_\a p(\a) \frac{\t^\a}{|\a|!} = \exp \biggl( \sum_{j=1}^\infty
\frac{t_j}{j!} \biggr) .
$$
\end{lemma}
\begin{proof}
Indeed, $p(\a)$ is the number of automorphisms of the set with $|\a|$
elements divided by the number of automorphisms of such a partition, namely
$$
p(\a) = |\a|! \Bigm/ \prod_{j=1}^\infty j!^{a_j}a_j! ,
$$
from which the lemma follows.
\end{proof}

Let $f$ be a power series
$$
f(t) = \sum_{k=1}^\infty \frac{f_kt^k}{k!} .
$$
Define the partial Bell polynomials $B_{n,k}$ by the generating function
$$
\exp(xf(t)) = \sum_{n=0}^\infty \frac{t^n}{n!} \sum_{k=0}^n
B_{n,k}(f_1,\dots,f_n) x^k .
$$
Setting $t_j=xt^jf_j$ in the generating function $B(\t)$ of Lemma
\ref{partitions}, we obtain the explicit formula
\begin{equation} \label{Bell}
B_{n,k}(f_1,\dots,f_n) = \sum_{J\in\S(n,k)} \prod_{i=1}^k f_{|J_i|} .
\end{equation}
(See Ex.\ 2.11 of Macdonald \cite{Macdonald}.) In particular, the partial
Bell polynomials have positive integral coefficients.

\begin{proposition} \label{inverse}
If $g$ is the inverse power series to $f$ (that is, $g(f(t))=t$), then
the matrices $F_{nk}=B_{n,k}(f_1,\dots,f_n)$ and
$G_{nk}=B_{n,k}(g_1,\dots,g_n)$ are inverse to each other.
\end{proposition}
\begin{proof}
The matrix $F$ is the transition matrix between the bases $(t^k/k!\mid
k\ge0)$ and $(f(t)^k/k!\mid k\ge0)$ of $\Q[t]$. Its inverse $F^{-1}$ is
thus the transition matrix between the bases $(f(t)^k/k!\mid k\ge0)$ and
$(t^k/k!\mid k\ge0)$. Changing variables from $t$ to $g(t)$, the result
follows.
\end{proof}

\subsection{Stirling numbers of the first kind}
The Stirling number of the first kind $s(n,k)$ may be defined as
$(-1)^{n-k}$ times the number of permutations on $n$ letters with $k$
cycles. A permutation of the set $\{1,\dots,n\}$ is the same thing as a
partition of $n$, together with a cyclic order on each part of the
partition. Since a set of cardinality $i$ has $(i-1)!$ cyclic orders, we
see that
\begin{equation} \label{link}
s(n,k) = \sum_{J\in\S(n,k)} \prod_{i=1}^n (-1)^{|J_i|-1} (|J_i|-1)! .
\end{equation}
Applying \eqref{Bell} with $f_k=(-1)^{k-1}(k-1)!$ (or $f(t)=\log(1+t)$), we
see that
\begin{equation} \label{first}
1 + \sum_{n=1}^\infty \sum_{k=1}^n s(n,k) \frac{t^nx^k}{n!} = (1+t)^x .
\end{equation}
In particular, for $n\ge1$,
\begin{equation} \label{descending}
\sum_{k=1}^n s(n,k) x^k = x(x-1)\dots(x-n+1) .
\end{equation}

\subsection{Stirling numbers of the second kind}
The number $S(n,k)$ of partitions of $n$ with $k$ parts (i.e.\ the
cardinality of $\S(n,k)$) is called a Stirling number of the second
kind. The special case of \eqref{Bell} with $f(t)=e^t-1$ (and hence $f_k=1$
for all $k$) shows that the Stirling numbers of the second kind have
generating function
\begin{equation} \label{second}
\sum_{n=1}^\infty \sum_{k=1}^n S(n,k) \frac{t^nx^k}{n!} = e^{x(e^t-1)} - 1
.
\end{equation}

For the reader's edification, we display the first few rows of the matrices
of first and second Stirling numbers:
$$
s = \begin{bmatrix}
1 & 0 & 0 & 0 & 0 & 0 \\
-1 & 1 & 0 & 0 & 0 & 0 \\
2 & -3 & 1 & 0 & 0 & 0 \\
-6 & 11 & -6 & 1 & 0 & 0 \\
24 & -50 & 35 & -10 & 1 & 0 \\
\hdotsfor{5} & \ddots
\end{bmatrix}
\quad\text{and}\quad
S = \begin{bmatrix}
1 & 0 & 0 & 0 & 0 & 0 \\
1 & 1 & 0 & 0 & 0 & 0 \\
1 & 3 & 1 & 0 & 0 & 0 \\
1 & 7 & 6 & 1 & 0 & 0 \\
1 & 15 & 25 & 10 & 1 & 0 \\
\hdotsfor{5} & \ddots
\end{bmatrix} .
$$

Applying Proposition \ref{inverse} to the functions $f(t)=e^t-1$ and
$g(t)=\log(1+t)$, we see that the matrices $s$ and $S$ formed from the
numbers $s(n,k)$ and $S(n,k)$ are inverse to each other:
\begin{equation} \label{Stirling}
\sum_{n=1}^\infty s(j,n) S(n,k) = \delta(j,k)
\end{equation}

We may rewrite \eqref{descending} in the form
$$
\sum_{j=1}^n s(j,n) x^{-n} = x^{-j} (1-x)\dots(1-(j-1)x) .
$$
{}From \eqref{Stirling}, it now follows easily that
$$
\sum_{n=k}^\infty S(n,k) x^n = \frac{x^k}{(1-x)(1-2x)\dots(1-kx)} ,
$$
or equivalently,
$$
S(n,k) = \sum_{1\le i_1\le\dots\le i_k\le n} i_1\dots i_k .
$$
This is also not difficult to prove directly, by induction on $n$.

\subsection{The cohomology of the configuration spaces $\FF(\C,n)$}
Let $H^\bull(\FF(\C,n),\Z)$ be the cohomology of the configuration space of
the complex line. Given distinct $j$ and $k$ in $\{1,\dots,n\}$, let
$\om_{jk}\in H^1(\FF(\C,n),\Z)$ be the integral cohomology class
represented by the closed differential form
$$
\Om_{jk} = \frac{1}{2\pi i} \frac{d(z_j-z_k)}{z_j-z_k} .
$$
By induction on $n$, Arnold shows in \cite{Arnold} that the cohomology ring
$H^\bull(\FF(\C,n),\Z)$ is generated by the classes $\om_{jk}$, subject to
the relations $\om_{jk}=\om_{kj}$ and
\begin{equation} \label{Arnold}
\om_{ij}\om_{jk} + \om_{jk}\om_{ki} + \om_{ki}\om_{ij} = 0 .
\end{equation}
The action of the group $\SS_n$ on the configuration space $\FF(\C,n)$
induces an action on the cohomology ring $H^\bull(\FF(\C,n),\Z)$, which
permutes the generators, by the formula
$$
\sigma\*\om_{ij} = \om_{\sigma(i)\sigma(j)} .
$$

Using the above presentation of $H^\bull(\FF(\C,n),\Z)$, Arnold shows that
$H^{n-k}(\FF(\C,n),\Z)$ is a free abelian group of rank
$(-1)^{n-k}s(n,k)$. This motivates the definition of a graded
$\SS_n$-module $\s(n,k)$, with
$$
\s(n,k)^i = \begin{cases} H^i(\FF(\C,n),\Z) , & i=n-k , \\
0 . & \text{otherwise,} \end{cases}
$$
We may think of $\s(n,k)$ as a lift of the Stirling number $s(n,k)$ to the
category of graded $\SS_n$-modules.

Denote by $\l(n)$ the graded $\SS_n$-module $\s(n,1)$. More generally, if
$\n$ is a finite set of cardinality $n$, let $\l(\n)$ be the graded
$\Aut(\n)$-module defined in the same way as $\l(n)$ but with the set
$\{1,\dots,n\}$ replaced by $\n$. It is isomorphic to $\l(n)$, but to
obtain an isomorphism, we must choose a total order on $\n$.

The following theorem is proved by Orlik and Solomon for general hyperplane
arrangements \cite{OS} (see Theorem 4.21 of Orlik \cite{Orlik}). See also
Lehrer-Solomon \cite{LS} for the special case which we consider.
\begin{theorem} \label{Lehrer-Solomon}
There is a natural decomposition
\begin{equation} \label{decompose}
\s(n,k) = \bigoplus_{J\in\S(n,k)} \s(n,J) ,
\end{equation}
together with natural isomorphisms
$$
\s(n,J) \cong \bigotimes_{i=1}^k \l(J_i) .
$$
\end{theorem}
\begin{proof}
The graded $\SS_n$-module $H^\bull(\FF(\C,n),\Z)$ is spanned by monomials
in the generators $\om_{ij}$. To such a monomial, we associate a forest (a
graph each of whose components is a tree) with vertices the set
$\{1,\dots,n\}$, and with an edge between vertices $i<j$ if and only if the
generator $\om_{ij}$ occurs in the monomial. Such a forest determines a
partition of the set $\{1,\dots,n\}$. Let $\s(n,J)$ be the span of the
monomials associated to the partition $J$.

Since all of the forests related by application of one of Arnold's
relations \eqref{Arnold} give rise to the same partition, we see that
$\s(n,J)$ is well-defined. If $J$ has $k$ parts, its associated forest has
$n-k$ edges, and hence $\s(n,J)$ is a subgroup of $H^{n-k}(\FF(\C,n),\Z)$.
In particular, if $J$ is the unique partition in $\S(n,1)$, we see that
$\s(n,J)$ is generated by all trees with $n$ labelled vertices, modulo the
Arnold relations. From this, it is easy to see that
$$
\s(n,J) \cong \bigotimes_{i=1}^k \s(J_i,1) .
\qed$$
\def\qed{}
\end{proof}

The characters of the $\SS_n$-modules $\l(n)$ have been calculated by
Hanlon \cite{Hanlon1} and Stanley \cite{Stanley}. From their formula, one
may calculate the characters of $\s(n,k)$ for all $k$.
\begin{lemma} \label{Hanlon-Stanley}
The equivariant Euler characteristic of the graded $\SS_n$-module
$\l(n)$, evaluated at $\sigma\in\SS_n$, is given by the formula
$$
\chi_\sigma(\l(n)) = \begin{cases} \displaystyle
- \frac{\mu(d)}{n}(-d)^{n/d}(n/d)! , & \text{if $\sigma$ has $n/d$
cycles of length $d$,} \\
0 , & \text{otherwise.}
\end{cases}$$
\end{lemma}

\subsection{A differential on $H^\bull(\FF(\C,n),\Z)$}
We now study the differential
$$
\p : H^\bull(\FF(\C,n),\Z) \to H^{\bull-1}(\FF(\C,n),\Z)
$$
of the algebra $H^\bull(\FF(\C,n),\Z)$ associated to the diagonal action of
the multiplicative group $\C^\times$ on $\FF(\C,n)$; it is given by capping
with the fundamental class of the circle $\U(1)\subset\C^\times$. It
follows from the definition of $\om_{ij}$ that $\p\om_{ij}=1$. One can
easily check that $\p$ is well-defined, by showing that the differential of
the relation \eqref{Arnold} vanishes:
$$
\p \bigl( \om_{ij}\om_{jk} + \om_{jk}\om_{ki} + \om_{ki}\om_{ij} \bigr)
=  \bigl( \om_{jk} - \om_{ij} \bigr) + \bigl( \om_{ki} - \om_{jk} \bigr) +
\bigl( \om_{ij} - \om_{ki} \bigr) = 0 .
$$
The following lemma reflects the fact that the action of $\C^\times$ on
$\FF(\C,n)$ is free if $n>1$, and that the resulting principal fibration is
trivial.
\begin{lemma} \label{acyclic}
If $n>1$, the complex $(H^\bull(\FF(\C,n),\Z),\p)$ is acyclic.
\end{lemma}
\begin{proof}
Let $H$ denote the operator of multiplication by $\om_{12}$. Since
$\p\om_{12}=1$, we see that $\p\*H+H\*\p$ equals the identity operator,
proving acyclicity.
\end{proof}

If $J\in\S(n,j)$ and $K\in\S(n,k)$ are partitions of $\{1,\dots,n\}$,
denote by $\p_{JK}$ the component of $\p$ mapping from $\s(n,J)$ to
$\s(n,K)$; thus, $\p_{JK}$ vanishes unless $k=j+1$.
\begin{lemma}
The differential $\p_{JK}$ vanishes unless $K\prec J$.
\end{lemma}
\begin{proof}
Let $\alpha$ be a monomial in the generators $\om_{ij}$ of
$H^\bull(\FF(\C,n),\Z)$. By the definition of $\p$, $\p\alpha$ is a sum of
terms, in each of which one of the factors $\om_{ij}$ occurring in $\alpha$
is omitted. Such a term corresponds to partition of $\{1,\dots,n\}$ in
which $i$ and $j$ are no longer equivalent: in other words, the partition
associated to the new monomial is a refinement of the partition associated
to $\alpha$.
\end{proof}

\section{Resolving sheaves on configuration spaces}

Before turning to the construction of resolutions of sheaves over
configuration spaces, we explain by an informal argument why one expects
Stirling numbers to arise in the construction.

Let $\pi:X\to S$ be a continuous map of locally compact topological spaces,
and let $X^n/S$ be the $n$th fibred power of $X$ with itself, defined
inductively by $X^0/S=S$ and
$$
X^{n+1}/S = (X^n/S) \times_S X .
$$
Denote by $\pi(n):X^n/S\to S$ the projection to $S$.

The space $X^n/S$ has a stratification, with strata indexed by the poset
of partitions $J$ of $\{1,\dots,n\}$: the stratum associated to a partition
$J$ is given by
$$
\FF(X/S,J) = \{ (x_1,\dots,x_n) \in X^n/S \mid \text{$x_i=x_j$ iff
$i\sim_Jj$} \} .
$$
A stratum $\FF(X/S,K)$ lies in the closure of $\FF(X/S,J)$ if and only if
$J\prec K$; the closure of $\FF(X/S,J)$ is the diagonal
$$
X^J/S = \{ (x_1,\dots,x_n) \in X^n/S \mid \text{$x_i=x_j$ if $i\sim_Jj$} \}
.
$$
If $J\in\S(n,k)$, denote by $i(J):X^J/S\hookrightarrow X^n/S$ the diagonal
embedding. If $\CF$ is a sheaf on $X^n/S$, denote by $\CF(J)$ the sheaf
$i(J)_!i(J)^*\CF$ on $X^n/S$.

If $J\in\S(n,k)$, $\FF(X/S,J)$ is isomorphic to $\FF(X/S,k)$; thus, we may
represent the above stratification of $X^n/S$ (in)formally as
$$
X^n/S = \coprod_{k=1}^n S(n,k) \* \FF(X/S,k) .
$$
Equation \eqref{Stirling} leads us to expect that there is a ``virtual
stratification'' of $\FF(X/S,n)$, of the form
\begin{equation} \label{haha}
\FF(X/S,n) = \coprod_{k=1}^n s(n,k) \* X^k/S .
\end{equation}
Rewritten in terms of generating functions, this becomes
$$
\sum_{n=0}^\infty \frac{t^n[\FF(X/S,n)]}{n!} = \sum_{k=0}^\infty
\frac{\log(1+t)^n[X^n/S]}{n!} = (1+t)^{[X/S]} ,
$$
where we think of the symbol $[X^n/S]$ as the $n$th power of $[X/S]$, as
indeed it is in the Grothendieck group of motivic sheaves on $S$.

We may make sense of \eqref{haha} using complexes of sheaves on
$X^n/S$. Let $j(n):\FF(X/S,n)\hookrightarrow X^n/S$ be the open embedding
of the configuration space in $X^n/S$. There is a natural resolution
$\Resolve^\bull(X/S,n,\CF)$ of $j(n)_!j(n)^*\CF$, whose underlying graded
sheaf has the form
\begin{equation} \label{resolve}
\Resolve^{n-k}(X/S,n,\CF) = \bigoplus_{J\in\S(n,k)} \Hom(\s(n,J),\CF(J)) .
\end{equation}
For example, if $n=2$, we recover \eqref{diagonal} while if $n=3$, we
obtain the resolution
$$
0 \rightarrow j(3)_!j(3)^*\CF \rightarrow \CF(1,2,3) \rightarrow \CF(12,3)
\oplus \CF(13,2) \oplus \CF(23,1) \rightarrow \CF(123)\oplus\CF(123)
\rightarrow 0 .
$$

In the special case that $\CF$ is a constant sheaf, we may interpret
$\CF(J)$ as a copy of the diagonal $X^J$, which is isomorphic to $X^k$ when
$J\in\S(n,k)$; replacing the $\SS_n$-module $\s(n,J)$ by its Euler
characteristic and bearing in mind \eqref{decompose}, or its numerical
version \eqref{link}, we are led to \eqref{haha}.

\subsection{Construction of the resolution}
If $J\prec K$ are partitions of $\{1,\dots,n\}$, denote by $i(J,K)$ the
inclusion $X^K/S\hookrightarrow X^J/S$, and by $i(J,K)^*:\CF(J)\to\CF(K)$
the induced map of sheaves. Let $\Resolve^\bull(X/S,n,\CF)$ be the complex
of sheaves \eqref{resolve}, with differential
\begin{equation} \label{differential}
d = \sum_{J\prec K} \p_{KJ}^* \o i(J,K)^* .
\end{equation}
For example, $\Resolve^0(X/S,n,\CF)\cong\CF$, while $\Resolve^1(X/S,n,\CF)$
is the direct sum
$$
\Resolve^1(X/S,n,\CF) = \bigoplus_{1\le k<l\le n}
\CF(kl,1,\dots,\widehat{k},\dots,\widehat{l},\dots,n) ,
$$
since $\dim\s(n,J)=1$ for all $J\in\S(n,n-1)$. In particular,
$j^*\Resolve^1(X/S,n,\CF)=0$.

Denote by $\eta:j(n)_!j(n)^*\To\Id$ the unit of the adjunction between
$j(n)_!$ and $j(n)^*$; it induces a map, also denoted by $\eta$, from
$j(n)_!j(n)^*\CF$ to $\CF=\Resolve^0(X/S,n,\CF)$. The composition of arrows
$$
j(n)_!j(n)^*\CF \xrightarrow{\eta} \Resolve^0(X/S,n,\CF) \xrightarrow{d}
\Resolve^1(X/S,n,\CF)
$$
is zero, showing that $\eta:j(n)_!j(n)^*\CF\to\Resolve^\bull(X/S,n,\CF)$ is
a morphism of complexes.
\begin{theorem} \label{Resolve}
The morphism $\eta:j(n)_!j(n)^*\CF\to\Resolve(X/S,n,\CF)$ is a
quasi-isomorphism.
\end{theorem}
\begin{proof}
We apply the following lemma.
\begin{lemma} \label{strata}
Let $X$ be a locally compact stratified space with strata $\{X_J\}$, and
let $j(J)$ be the locally closed embedding of the stratum $X_J$ in
$X$. Then a map of complexes of sheaves $\eta:\CF_1\to\CF_2$ on $X$ is a
quasi-isomorphism if and only if $\eta : j(J)_!j(J)^*\CF_1 \to j(J)_!j(J)^*
\CF_2$ is a quasi-isomorphism for all strata.
\end{lemma}
\begin{proof}
If there is only one stratum, the lemma is a tautology. We now argue by
induction on the number of strata. Let $X_J$ be an open stratum of $X$, and
let $Z$ be its complement in $X$, with closed embedding
$i:Z\hookrightarrow X$. Consider the diagram
$$\begin{diagram}
\node{0} \arrow{e} \node{j(J)_!j(J)^*\CF_1} \arrow{e} \arrow{s,l}{\eta}
\node{\CF_1} \arrow{e} \arrow{s,l}{\eta} \node{i_!i^*\CF_1} \arrow{e}
\arrow{s,l}{\eta} \node{0} \\
\node{0} \arrow{e} \node{j(J)_!j(J)^*\CF_2} \arrow{e} \node{\CF_1} \arrow{e}
\node{i_!i^*\CF_1} \arrow{e} \node{0}
\end{diagram}$$
Since the rows are exact, we conclude by the five-lemma that
$\eta:\CF_1\to\CF_2$ is a quasi-isomorphism if and only if
$\eta:j(J)_!j(J)^*\CF_1\to j(J)_!j(J)^*\CF_2$ and $\eta:i_!i^*\CF_1\to
i_!i^*\CF_2$ are. By the induction hypothesis, $\eta:i_!i^*\CF_1\to
i_!i^*\CF_2$ is a quasi-isomorphism if and only if
$\eta:j(K)_!j(K)^*\CF_1\to j(K)_!j(K)^*\CF_2$ are for all $K\ne J$; this
proves the induction step.
\end{proof}

If $J$ is a partition of $\{1,\dots,n\}$, let
$j(J):\FF(X/S,J)\hookrightarrow X^n/S$ be the inclusion of the locally
closed subspace $\FF(X/S,J)$. By the base change theorem,
$$
j(J)_!j(J)^*j(n)_!j(n)^*\CF \cong \begin{cases} j(n)_!j(n)^*\CF , &
\text{if $J$ is the unique partition in $\S(n,n)$,} \\
0 , & \text{otherwise.}
\end{cases}$$
The base change theorem we apply here is (1.2.1) of Verdier
\cite{Verdier2}. (In the language of Section 3 of this paper, it says that
the natural transformation $\phi:g^*t_!\To s_!f^*$ of Definition
\ref{MACKEY} (ii) is an isomorphism.)

Let the parts of $J$ be $\{J_1,\dots,J_k\}$, in no particular order, and
let $n_i$ be the cardinality of $J_i$. Applying Theorem
\ref{Lehrer-Solomon}, we see that
\begin{align*}
j(J)_!j(J)^*\Resolve^\bull(X/S,n,\CF) & \cong \Hom \biggl(
\bigoplus_{K\prec J} \s(n,K) , j(J)_!j(J)^*\CF \biggr) \\
& \cong \Hom \biggl( \bigotimes_{i=1}^k H^{n_i-\bull}(\FF(\C,n_i),\Z) ,
j(J)_!j(J)^*\CF \biggr) .
\end{align*}
The differential on this complex of sheaves is induced by the differentials
on the factors $H^{n_i-\bull}(\FF(\C,n_i),\Z)$, and hence by Lemma
\ref{acyclic} is acyclic if $n_i>1$ for any $i$.

We see that the hypotheses of Lemma \ref{strata} are fulfilled: if $J$ is a
partition of $\{1,\dots,n\}$, $\eta:j(J)_!j(J)^*\CF\to
j(J)_!j(J)^*\Resolve(X/S,n,\CF)$ is a quasi-isomorphism, since the two
complexes are equal if $J\in\S(n,n)$, while they are both acyclic otherwise.
\end{proof}

\section{Mackey $2$-functors}

In this section, we axiomatize those properties of the $2$-functor
associating to a variety its derived category of mixed Hodge modules which
will be used in constructing the analogue of the resolution
$\Resolve^\bull(X/S,n,\CF)$ for mixed Hodge modules. It turns out that we
need a natural analogue for $2$-functors of Dress's Mackey functors
\cite{Dress}.

Impatient readers may skip to Section \ref{exact}: all they need to know
about the Mackey $2$-functor underlying the theory of mixed Hodge modules
is that the usual properties of the functors $f_!$ and $f^*$ for locally
closed \embeddings hold, such as the base change theorem (in particular, we
make no use of Verdier duality). In Section \ref{exact}, we impose
sufficient additional hypotheses on these functors to allow us to construct
\Cech-type resolutions of $j_!j^*\CF$ when $j$ is an open \embedding and of
$i_!i^*\CF$ when $i$ is a closed \embedding. The analogue of
$\Resolve(X/S,n,\CF)$ for mixed Hodge modules is defined by replacing the
sheaf $i(J)_!i(J)^*\CF$ by this \Cech complex.

In all $2$-categories which occur in this paper, the $2$-morphisms are
invertible. We often consider $2$-functors which have a category as their
domain, which is thought of as a $2$-category all of whose $2$-morphisms
are identities.

\subsection{Mackey functors}
Mackey functors were introduced by Dress \cite{Dress} as an axiomatization
of induction in the theory of group representations. The motivating example
is the functor $G\mapsto R(G)$ on the category of finite groups, which
assigns to a group $G$ its virtual representation ring. Given a morphism
$f:G\to H$ of finite groups, there is a contravariant map $f^\bull:R(H)\to
R(G)$, pull-back along $f$, and a covariant map $f_\bull:R(G)\to R(H)$
generalizing induction:
$$
f_\bull V = ( \C[G] \o V )^H .
$$
These functors satisfy the Mackey double coset formula, which says that
given a Cartesian square of finite groups
$$\begin{diagram}
\node{G_1} \arrow{e,t}{f} \arrow{s,l}{s} \node{G_2} \arrow{s,r}{t} \\
\node{H_1} \arrow{e,t}{g} \node{H_2}
\end{diagram}$$
we have the equality $g^\bull t_\bull=s_\bull f^\bull$.

\subsection{Mackey $2$-functors}
The group of virtual representations $R(G)$ is the Gro\-then\-dieck group
of the category $\Proj(G)$ of finite-dimensional representations. The
$2$-functor $G\mapsto\Proj(G)$ satisfies axioms which are the natural
analogue for $2$-functors of the notion of a Mackey functor. Broadly
speaking, a Mackey $2$-functor is a $2$-functor satisfying these axioms,
along with an additivity axiom whose analogue for Mackey functors is only
meaningful if we allow $G$ to be a groupoid.

The concept of a Mackey $2$-functor is not new: it was introduced by
Deligne in Expos\'e XVII of SGA~4 \cite{SGA4} (though not under this name).

\begin{definition} \label{MACKEY}
A Mackey $2$-functor from a category $\Cat$ to a $2$-category $\TT$
consists of a pair of $2$-functors $\D^\bull:\Cat^\op\to\TT$ (where
$\Cat^\op$ is the opposite of $\Cat$) and $\D_\bull:\Cat\to\TT$, such that
\begin{enumerate}
\item if $X$ is an object of $\Cat$, the objects $\D^\bull(X)$ and
$\D_\bull(X)$ are identical --- we denote this object by $\D(X)$, and if
$f:X\to Y$ is a morphism of $\Cat$, we denote the $1$-morphism
$\D^\bull(f):\D(Y)\to\D(X)$ by $f^\bull$ and the $1$-morphism
$\D_\bull(f):\D(X)\to\D(Y)$ by $f_\bull$;
\item (base change) to each Cartesian square
$$\begin{diagram}
\node{X_1} \arrow{e,t}{f} \arrow{s,l}{s} \node{X_2} \arrow{s,r}{t} \\
\node{Y_1} \arrow{e,t}{g} \node{Y_2}
\end{diagram}$$
in $\Cat$ is associated a natural $2$-morphism $\phi:g^\bull t_\bull\To
s_\bull f^\bull$, such that given a diagram each square of which is
Cartesian
$$\begin{diagram}
\node{X_1} \arrow{e,t}{f} \arrow{s,l}{s} \node{X_2} \arrow{e,t}{f'}
\arrow{s,l}{t} \node{X_3} \arrow{s,l}{u} \\
\node{Y_1} \arrow{e,t}{g} \arrow{s,l}{s'} \node{Y_2} \arrow{e,t}{g'}
\arrow{s,l}{t'} \node{Y_3} \\
\node{Z_1} \arrow{e,t}{h} \node{Z_2}
\end{diagram}$$
the $2$-morphism $\phi$ associated to the top (resp.\ left) pair of squares
is the $2$-com\-po\-si\-ti\-on of the $2$-morphisms associated to the squares
from which it is formed;
\item (additivity) there are $2$-morphisms of bifunctors
$\alpha:\D^\bull(X\coprod Y)\To \D^\bull(X)\oplus \D^\bull(Y)$ and
$\beta:\D_\bull(X\coprod Y)\To \D_\bull(X)\oplus \D_\bull(Y)$, such that
$$
\alpha_{X,Y} = \beta_{X,Y} : {\textstyle\D(X\coprod Y)} \To \D(X)\oplus \D(Y) .
$$
\end{enumerate}
\end{definition}

\subsection{Exact Mackey $2$-functors} \label{exact}
Let $\Var$ be the category of quasi-projective complex varieties, with
morphisms the locally closed \embeddings. Let $\TT$ be the $2$-category
whose objects are $t$-categories (Beilinson-Bernstein-Deligne \cite{BBD}),
whose $1$-morphisms are right $t$-exact functors possessing a right
adjoint, and whose $2$-morphisms are natural isomorphisms.
\begin{definition}
An exact Mackey $2$-functor on $\Var$ is a Mackey $2$-functor with values in
$\TT$ such that
\begin{enumerate}
\item for closed \embeddings $i:Z\hookrightarrow X$, $i^\bull$ has right
adjoint $i_\bull$, and $i_\bull$ is fully faithful;
\item for open \embeddings $j:U\hookrightarrow X$, $j_\bull$ 
has right adjoint $j^\bull$ and is fully faithful;
\item if $i:Z\hookrightarrow X$ is a closed \embedding, and
$j:U\hookrightarrow X$ is the open \embedding of the complement
$U=X\setminus Z$, there is an exact triangle $(j_\bull j^\bull V,V,i_\bull
i^\bull V)$ for each object $V$ of $\D(X)$.
\item if $j$ is an affine open \embedding, the functor $j_\bull$ is
$t$-exact.
\end{enumerate}
\end{definition}

We have in mind three examples of exact Mackey $2$-functors.
\begin{example}
The $2$-functors $X\longmapsto(D^b(X),f_\bull=f_!,f^\bull=f^*)$ assigning to
$X$ the derived category of sheaves of abelian groups with the usual
$t$-structure.
\end{example}

\begin{example}
The $2$-functors $X\longmapsto({}^p\!D^b_c(X),f_\bull=f_!,f^\bull=f^*)$
assigning to $X$ the derived category of bounded complexes of sheaves with
constructible cohomology, together with the $t$-structure associated to the
\emph{middle} perversity $p$ (axiom (iv) is Corollaire 4.1.3 of
\cite{BBD}).
\end{example}

\begin{example}
The $2$-functors $X\longmapsto(D^b(\MHM(X)),f_\bull=f_!,f^\bull=f^*)$
assigning to $X$ the derived category of mixed Hodge modules with the
natural $t$-structure.
\end{example}

The Grothendieck group $K_0(\DD)$ of a triangulated category $\DD$ is the
abelian group generated by the isomorphism classes of objects of $\DD$
(where we assume that these form a set); we impose the relation
$[V]\sim[U]+[W]$ for all exact triangles $(U,V,W)$ in $\D$. For example, if
$\DD$ is the derived category $D^b(\Ab)$ of bounded complexes in an abelian
category $\Ab$, then $K_0(\DD)$ may be identified with the Grothendieck
group $K_0(\Ab)$ of the abelian category $\Ab$, which is the abelian group
generated by the isomorphism classes of objects of $\Ab$, with the relation
$[V]\sim[U]+[W]$ whenever $V$ is an extension of $W$ by $U$.

Although we will not need it, the following result goes some way towards
justifying our introduction of Mackey $2$-functors.
\begin{proposition}
Let $\D$ be an exact Mackey $2$-functor, and let $\K$ be the composition of
$\D$ with the functor $K_0$ from $\TT$ to $\Ab$. Then $\K$ is a Mackey
functor on $\Var$.
\end{proposition}

\subsection{\Cech complexes}
Let $j:U\hookrightarrow X$ be an open \embedding and let $\CU=\{U_i\}_{1\le
i\le d}$ be a finite cover of $U$ by affine open \embeddings
$j_i:U_i\hookrightarrow X$. For example, we might take the $U_i$ to be
complements of Cartier divisors in $X$. Using these data, we now define a
\Cech-type resolution of $j_\bull j^\bull \CF$, where $\D$ is an exact
Mackey $2$-functor. (See also Proposition 2.19 of Saito \cite{Saito:mixed}
and Section 3.4 of Beilinson \cite{Beilinson}.)
\begin{definition}
The \Cech-complex $\CC_\bull(X,\CU,\CF)$ is the graded object of $\D(X)$
$$
\CC_k(X,\CU,\CF) = \bigoplus_{i_0<\dots<i_k} (j_{i_0\dots i_k})_\bull
(j_{i_0\dots i_k})^\bull \CF ,
$$
where $j_{i_0\dots i_k}$ is the open \embedding of $U_{i_0\dots i_k} =
\bigcap_{\ell=0}^k U_i$ in $X$. Its differential is the sum of maps
$$
\p = \sum_{\ell=0}^k (-1)^\ell \p_\ell : \CC_k(X,\CU,\CF) \to
\CC_{k-1}(X,\CU,\CF) .
$$
Here, $\p_\ell : (j_{i_0\dots i_k})_\bull(j_{i_0\dots i_k})^\bull\CF \to
(j_{i_0\dots\widehat{\imath_\ell}\dots i_k})_\bull
(j_{i_0\dots\widehat{\imath_\ell}\dots i_k})^\bull\CF$ is induced by the
adjunction $q_\bull q^\bull\To\Id$ associated to the open \embedding
$q:U_{i_0\dots i_k}\hookrightarrow U_{i_0\dots\widehat{\imath_\ell}\dots
i_k}$.
\end{definition}

If $\DD$ is a $t$-category, let $H^0(\DD)$ be its heart. Recall from
Section 3.1.9 of \cite{BBD} the realization functor
$$
\real : D^b(H^0(\DD)) \to \DD ;
$$
this is an exact functor mapping bounded complexes $\CC_\bull(H^0(\DD))$ to
objects $\real(\CC_\bull)$ in $\DD$.

\begin{proposition} \label{Cech}
Let $\D$ be an exact Mackey $2$-functor, let $j:U\hookrightarrow X$ be an
open \embedding, let $i:Z\hookrightarrow X$ be the closed \embedding of the
complement $Z=X\setminus U$, and let $\CU=\{U_i\}_{0\le i\le d}$ be a cover
of $U$ by affine open \embeddings $j_i:U_i\hookrightarrow X$. Then
\begin{enumerate}
\item $\real\bigl(\CC_\bull(X,\CU,\CF)\bigr)\in\Ob\D(X)^{[-d,0]}$ is
isomorphic to $j_\bull j^\bull \CF$;
\item
$\real\bigl(\Cone\bigl(\CC_\bull(X,\CU,\CF)\to\CF\bigr)\bigr) \in
\Ob\D(X)^{[-d-1,0]}$ is isomorphic to $i_\bull i^\bull \CF$.
\end{enumerate}
\end{proposition}
\begin{proof}
Note that (ii) is implied by (i) and the five-lemma: in the diagram
$$\begin{diagram}
\divide\dgARROWLENGTH by 2
\node{0} \arrow[1]{e} \node[1]{\CF} \arrow[2]{e}
\arrow[2]{s,l}{\simeq} \node[2]{i_\bull i^\bull\CF} \arrow[2]{e}
\arrow[2]{s} \node[2]{j_\bull j^\bull\CF[1]} \arrow[1]{e}
\arrow[2]{s,l}{\simeq} \node[1]{0} \\[2]
\node{0} \arrow[1]{e} \node[1]{\CF}
\arrow[2]{e} \node[2]{\Cone\bigl(\CC_\bull(X,\CU,\CF)\to\CF\bigr)}
\arrow[2]{e} \node[2]{\CC_\bull(X,\CU,\CF)[1]} \arrow[1]{e}
\node[1]{0}
\end{diagram}$$
the top row is exact by axiom (iii) for exact Mackey $2$-functors, and the
bottom row is obviously exact.

We now prove (i) by induction on $d$: for $d=0$, $\CC_\bull(X,\CU,\CF)\cong
j_\bull j^\bull\CF$, and the proposition is a tautology.

The open subset $U^0=\bigcup_{i=1}^dU_i$ of $X$ has cover
$\CU^0=\{U_i\}_{1\le i\le d}$. Let $j^0$ and $j_0$ be the open \embeddings
of $U^0$ and $U_0$ in $X$, and define $\CF^0=(j^0)_\bull(j^0)^\bull\CF$ and
$\CF_0=(j_0)_\bull(j_0)^\bull\CF$. Let $p$ be the locally closed \embedding
of $U\setminus U^0=U_0\setminus U^0$ in $X$. We now form the diagram
$$\begin{diagram}
\divide\dgARROWLENGTH by 4
\node{0} \arrow[2]{e} \node[2]{\CF^0} \arrow[2]{e}
\arrow[3]{s,l}{\simeq} \node[2]{j_\bull j^\bull\CF} \arrow[3]{e}
\arrow[3]{s} \node[3]{p_\bull p^\bull\CF} \arrow[3]{e}
\arrow[3]{s,l}{\simeq} \node[3]{0} \\[3]
\node{0} \arrow[2]{e} \node[2]{\CC_\bull(X,\CU^0,\CF^0)} \arrow[2]{e}
\node[2]{\CC_\bull(X,\CU,\CF)} \arrow[3]{e}
\node[3]{\textstyle\Cone(\CC_\bull(X,\CU^0\cap U_0,\CF_0)\to\CF_0)}
\arrow[3]{e} \node[3]{0}
\end{diagram}$$
The lower row is defined by dividing the summands $(j_{i_0\dots i_k})_\bull
(j_{i_0\dots i_k})^\bull\CF$ of $\CC(X,\CU,\CF)$ into two classes:
\begin{enumerate}
\item if $i_0>0$, the term $(j_{i_0\dots i_k})_\bull (j_{i_0\dots
i_k})^\bull\CF$ is a summand of
$\CC_\bull(\CU^0,(j^0)_\bull(j^0)^\bull\CF)$;
\item if $i_0=0$, the term $(j_{0i_1\dots i_k})_\bull (j_{0i_1\dots
i_k})^\bull\CF$ is a summand of $\Cone_k(\CC_\bull(\CU^0\cap
U_0,\CF_0)\to\CF_0)$.
\end{enumerate}
In particular, the bottom row is exact.

The top row is exact by axiom (iii) for exact Mackey $2$-functors, applied
to the closed \embedding of $U\setminus U^0$ in $U$, while the outer
vertical arrows are quasi-isomorphisms by the induction hypothesis. The
proposition now follows by the five-lemma.
\end{proof}

\subsection{The resolution $\Resolve_\CU^\bull(X/S,n,\CF)$ for mixed Hodge
modules} Using the \Cech-com\-plex\-es $\CC_\bull(X,\CU,\CF)$, we now
construct a resolution of $j_!j^*\CF$, where $\CF$ is a mixed Hodge
module. The functor $i^*$ is not $t$-exact on mixed Hodge modules for
general closed \embeddings $i$; for this reason, our construction depends
on the choice of an auxiliary cover $\CU$ of $\FF(X/S,2)$ by affine open
\embeddings $j_i:U_i\hookrightarrow X$. (If $X/S$ is a smooth family of
curves, we may take the cover to have one element $\CU=\{\FF(X/S,2)\}$,
since in that case, the diagonal in $X^2/S$ is a Cartier divisor.) We will
actually construct the resolution in the more general setting of exact
Mackey $2$-functors.

For $k,l\in\{1,\dots,n\}$ with $k\ne l$, let $\pi_{kl}:X^n/S\to X^2/S$ be
the morphism which projects onto the $k$th and $l$th factors. Define a
cover $\CU(J)$ of the complement of the diagonal $i(J):X^J/S\hookrightarrow
X^n/S$ by
$$
\CU(J) = \bigl\{ \pi_{kl}^{-1}(U_i) \mid \text{$k\sim_Jl$ and $U_i\in\CU$}
\bigr\} .
$$
The open \embeddings $\pi_{kl}^{-1}(U_i)\hookrightarrow X^n/S$ are affine,
since affine open \embeddings are preserved under base change (EGA II,
1.6.2 \cite{EGA}).

By Proposition \ref{Cech}, the realization of the complex
$\Cone(\CC_\bull(X^n/S,\CU(J),\CF)\to\CF)$ is quasi-isomorphic to
$i(J)_\bull i(J)^\bull \CF$. If $J\prec K$, the morphism
$$
i(J,K)^* : i(J)_\bull i(J)^\bull\CF \to i(K)_\bull i(K)^\bull\CF
$$
is induced by an inclusion of complexes
$$
i(J,K)^* : \CC_\bull(X^n/S,\CU(J),\CF) \to \CC_\bull(X^n/S,\CU(K),\CF) ,
$$
which exists because the cover $\CU(K)$ contains the open cover $\CU(J)$.

As in the case of sheaves, our resolution of $j_\bull j^\bull\CF$ is a sum
over partitions \eqref{resolve}; unlike that case, the result is a double
complex, and we must take the realization of its total complex to obtain an
object of $\D(X^n/S)$. Let
$$
\Resolve_\CU^{n-k,-j}(X/S,n,\CF) = \bigoplus_{J\in\S(n,k)}
\Hom(\s(n,J),\Cone_j(\CC_\bull(X^n/S,\CU(J),\CF)\to\CF)) .
$$
There are two differentials: the analogue of differential
\eqref{differential},
$$
d = \sum_{J\prec K} \p_{KJ}^* \o i(J,K)^* ,
$$
and the \Cech-differential $\Resolve_\CU^{n-k,-j}(X/S,n,\CF) \to
\Resolve_\CU^{n-k,1-j}(X/S,n,\CF)$.

We may identify $\Resolve_\CU^{0,\bull}(X/S,n,\CF)$ with $\CF$; thus, there
is a natural coaugmentation
$$
\eta : j_\bull j^\bull\CF \to \Resolve_\CU^{0,\bull}(X/S,n,\CF) .
$$

If $\sigma$ is a permutation, the cover $\CU(J)$ is carried into
$\CU(\sigma\*J)$ by the action of $\sigma$, so that $\sigma$ maps
$\CC_\bull(X^n/S,\CU(J),\CF)$ isomorphically to
$\CC_\bull(X^n/S,\CU(\sigma\*J),\CF)$. The differential of
$\Resolve_\CU^{\bull,\bull}(X/S,n,\CF)$ is invariant under the action of
$\sigma$, showing that $\Resolve_\CU^{\bull,\bull}(X/S,n,\CF)$ carries an
action of $\SS_n$. It is clear that the coaugmentation $\eta$ is
$\SS_n$-equivariant.

We now come to the main result of this paper; we omit the proof, since it
is essentially identical to that of Theorem \ref{Resolve}.
\begin{theorem} \label{RESOLVE}
Let $\D$ be an exact Mackey $2$-functor. If $\pi:X\to S$ is a
quasi-projective morphism of comples varieties, $\CF$ is an object of
$\D(X)$ and $\CU$ is a cover of $\FF(X/S,2)$ by affine open \embeddings,
the coaugmentation
$$
\eta:j_\bull j^\bull\CF \to \Resolve_\CU^{\bull,\bull}(X/S,n,\CF)
$$
induces an $\SS_n$-equivariant quasi-isomorphism between the objects
$\real\bigl(\Tot\Resolve_\CU^{\bull,\bull}(X/S,n,\CF)\bigr)$ and $j_\bull
j^\bull\CF$ in $\D(X^n/S)$.
\end{theorem}

\section{The $\SS_n$-equivariant relative Serre \polynomial of $j(n)_\bull
j(n)^\bull\CE^{\boxtimes n}$}

\subsection{The relative Serre \polynomial}
Let $\D$ be an exact Mackey $2$-functor, and let $\pi:X\to S$ be a morphism
of quasi-projective complex varieties. If $\CF$ is an object of $\D(X)$,
the \emph{relative Serre \polynomial} $\Serre_S(X,\CF)$ of $\CF$ is the
class of $\pi_\bull\CF$ in $\K(S)$, where $\K(S)$ is the Grothendieck group
$K_0(\D(S))$ of the triangulated category $\D(S)$.

This terminology is motivated by the special case in which $\D=D^b(\MHM)$
and $S=\Spec(\C)$. If we apply to $\Serre(X,\CF)$ the homomorphism
$\eps:K_0\bigl(\MHM(\Spec(\C))\bigr)\to\Z[t]$ defined by
$$
\eps(V) = \sum_{i,k} (-1)^i t^k \dim \gr^W_kV^i ,
$$
we obtain the Serre \polynomial of $\CF$,
$$
\eps\bigl(\Serre(X,\CF)\bigr) = \sum_{i,k} (-1)^i t^k \dim \gr^W_k
H^i_c(X,\CF) .
$$

Now suppose a finite group $\Gamma$ acts on $X$ and $Y$, and the morphism
$\pi$ is $\Gamma$-equivariant. If $\CF$ is an object of $\D^\Gamma(X)$, the
equivariant relative Serre \polynomial $\Serre_S^\Gamma(X,\CF)$ is the class
of $\pi_\bull\CF$ in $\K^\Gamma(S)$, where $\K^\Gamma(X)$ is the
Grothendieck group of $\Gamma$-equivariant objects in $\D(X)$.

When $\D$ is defined over a field $\k$ of characteristic $0$ and the
categories $\D(X)$ have tensor products, the associated Grothendieck groups
$\K^\Gamma(X)$ are $\lambda$-rings, by the arguments of \cite{I}. If in
addition the action of $\Gamma$ on $X$ is trivial, the Peter-Weyl theorem
(see Theorem 3.2 of \cite{I}) implies the isomorphism of $\lambda$-rings
$$
\K^\Gamma(X) \cong R(\Gamma) \o \K(X) ,
$$
where $R(\Gamma)$ is the virtual representation ring of $\Gamma$ (the
Grothendieck group of fin\-ite-di\-men\-sion\-al $\k[\Gamma]$-modules).

In this section, we calculate
$\Serre_S^{\SS_n}(\FF(X/S,n),j(n)^\bull\CE^{\boxtimes n})$. We obtain a
generalization of Theorem 5.6 of \cite{I}, which is the special case of
Theorem \ref{push} where $\D=\MHM$, $S=\point$ and $\CE=\1$ is the unit of
$\D(X)$.

The reader interested only in the case of mixed Hodge modules may skip the
next paragraph, whose r\^ole is to axiomatize the projection axiom.

\subsection{Green $2$-functors}
If $G$ is a finite group, $R(G)$ is not only an abelian group, but also a
commutative ring: furthermore, the functor $f^\bull$ preserves this
product, and we have the projection axiom, which says that if $f:G\to H$ is
a morphism of finite groups, there is a commutative diagram
$$\begin{diagram}
\node[3]{R(G)\o R(H)} \arrow{wsw,t}{1\o f^\bull} \arrow{ese,t}{f_\bull\o1} \\
\node{R(G)\o R(G)} \arrow{s} \node[4]{R(H)\o R(H)} \arrow{s} \\
\node{R(G)} \arrow[4]{e,b}{f_\bull} \node[4]{R(H)}
\end{diagram}$$
where the vertical arrows are multiplication in $R(G)$ and $R(H)$. Dress
calls a Mackey functor with these additional structures a Green functor.

The following definition gives the analogous condition for exact Mackey
$2$-functors.

\begin{definition}
An exact Green $2$-functor is an exact Mackey $2$-functor $\D$ on $\Var$
such that given a morphism $f:X\to Y$ in $\Var$, there is a natural
$2$-morphism
$$\divide\dgARROWLENGTH by2
\begin{diagram}
\node[3]{\D(X\times Y)} \arrow{wsw,t}{(1\times f)^\bull}
\arrow{ese,t}{(f\times1)_\bull} \\
\node{\D(X\times X)} \arrow{s,l}{i_X^\bull} \node[2]{\stackrel{\psi_f}\To}
\node[2]{\D(Y\times Y)} \arrow{s,r}{i_Y^\bull} \\
\node{\D(X)} \arrow[4]{e,b}{f_\bull} \node[4]{\D(Y)}
\end{diagram}$$
where $i_X:X\to X\times X$ and $i_Y:Y\to Y\times Y$ are the diagonal
\embeddings.

These $2$-morphisms must satisfy the condition that for any pair of
composable arrows $X\xrightarrow{f}Y\xrightarrow{g}Z$, the $2$-morphism
$\psi_{gf}$ equals the $2$-pasting of the following diagram:
$$
\divide\dgARROWLENGTH by5
\multiply\dgARROWLENGTH by3
\begin{diagram}
\node[5]{\D(X\times Z)} \arrow{wsw,t}{(1\times g)^\bull}
\arrow{ese,t}{(f\times1)_\bull} \\
\node[3]{\D(X\times Y)} \arrow{wsw,t}{(1\times f)^\bull}
\arrow{ese,t}{(f\times1)_\bull}
\node[2]{\overset{\phi}{\To}}
\node[2]{\D(Y\times Z)} \arrow{wsw,t}{(1\times g)^\bull}
\arrow{ese,t}{(g\times1)_\bull} \\
\node{\D(X\times X)} \arrow{s,l}{i_X^\bull} \node[2]{\stackrel{\psi_f}\To}
\node[2]{\D(Y\times Y)} \arrow{s,l}{i_Y^\bull} \node[2]{\stackrel{\psi_g}\To}
\node[2]{\D(Z\times Z)} \arrow{s,r}{i_Z^\bull} \\
\node{\D(X)} \arrow[4]{e,b}{f_\bull} \node[4]{\D(Y)}
\arrow[4]{e,b}{g_\bull} \node[4]{\D(Z)}
\end{diagram}$$
\end{definition}

\bigskip

All three examples of exact Mackey $2$-functors which we gave are exact
Green $2$-functors.

\subsection{A formula for
$\Serre_S^{\SS_n}\bigl(\FF(X/S,n),j(n)^\bull\CE^{\boxtimes n}\bigr)$} Let
$\D$ be an exact Green $2$-functor defined over $\Q$ (that is, all of the
$t$-categories $\D(X)$ are defined over $\Q$). If $\pi:X\to S$ is a
quasi-projective morphism of complex varieties and $\CE$ is an object of
$\D(X)$, we define
$$
\CE^{\boxtimes n} = \pi_1^\bull\CE \o \dots \o \pi_n^\bull\CE \in
\Ob\D^{\SS_n}(X^n/S) ,
$$
where $\pi_i:X^n/S\to X$ is the $i$th projection.

In calculating $\Serre_S^{\SS_n}(\FF(X/S,n),j(n)^\bull\CE^{\boxtimes n})$,
we make free use of the results of \cite{I}. Recall from loc.\ cit.\ that
if $R$ is a complete $\lambda$-ring, with decreasing filtration $F_iR$,
$i\ge0$, there is an operation
$$
\Exp(x) = \sum_{n=0}^\infty \sigma_n(x) = \exp\Bigl( \sum_{n=1}^\infty
\frac{1}{n} \psi_n(x) \Bigr) : F_1R \to 1+F_1R
$$
an analogue of the exponential, with inverse
$$
\Log(x) = \sum_{n=0}^\infty \frac{\mu(n)}{n} \log(\psi_n(x)) : 1 + F_1R \to
F_1R .
$$
Here, $\sigma_n$ is the $n$th $\sigma$-operation on $R$, and $\psi_n$ is
the $n$th Adams operation. For example, if $R=\Z\[q\]$ has its standard
$\lambda$-ring structure, then
$$
\Exp(a_1q+a_2q^2+a_3q^3+\dots) = (1-q)^{-a_1}(1-q^2)^{-a_2}(1-q^3)^{-a_3}
\dots.
$$

The ring of symmetric functions $\Lambda$ is the graded ring
$$
\Lambda = \varprojlim \Z[x_1,\dots,x_k]^{\SS_k} ,
$$
where the variables $x_i$ are assigned degree $1$. It is a polynomial ring
in the complete symmetric functions
$$
h_n = \sum_{i_1\le\dots\le i_n} x_{i_1}\dots x_{i_n} ,
$$
and may thus be identified with the free $\lambda$-ring on one generator
$h_1$, such that $\sigma_n(h_1)=h_n$ is the $n$th complete symmetric
function, and
$$
\psi_n(h_1) = p_n = \sum_i x_i^n
$$
is the $n$th power sum. Denote the abelian group of symmetric polynomials
of degree $n$ by $\Lambda_n$: it is free of rank $p(n)$, the number of
partitions of $n$.

If $R$ is a $\lambda$-ring, denote by $\Lambda\ohat R$ the complete tensor
product of $\Lambda$ with $R$, with filtration
$$
F_i(\Lambda\ohat R) = \Bigl\{ \sum_{n=i}^\infty a_i\o r_i \Big|
a_n\in\Lambda_n , r_n\in R \Bigr\} , \quad i\ge0 .
$$
(In particular, $\Lambda\ohat R=F_0(\Lambda\ohat R)$.)

If $V$ is an $\SS_n$-module defined over $\Q$, denote by $\ch(V)\in\Lambda$
its Frobenius characteristic; this is the degree $n$ symmetric function
given by the explicit expression
$$
\ch(V) = \frac{1}{n!} \sum_{\sigma\in\SS_n} \Tr_V(\sigma) p_\sigma ,
$$
where $p_\sigma$ is the monomial in the power sums obtained by taking one
factor $p_k$ for each cycle of $\sigma$ of length $k$. By the Peter-Weyl
theorem, this definition may be extended to $\SS_n$-modules in a $1$-ring
$\CR$ defined over $\Q$: there is a naturally defined isomorphism $\ch_n$
between the Grothendieck group of $\SS_n$-modules in $\CR$ and $\Lambda_n\o
K_0(\CR)$ (see Theorem 4.8 of \cite{I}).

\begin{theorem} \label{push}
Let $\D$ be an exact Green $2$-functor, and let $\pi:X\to S$ be a morphism
of quasi-projective complex varieties. If $\CE$ is an object of $\D(X)$,
the following equality holds in $\Lambda\ohat\K(S)$:
$$
\sum_{n=0}^\infty \Serre_S^{\SS_n}(\FF(X/S,n),j(n)^\bull\CE^{\boxtimes n})
= \Exp\biggl( \sum_{n=1}^\infty \frac{\mu(n)}{n} \Serre\bigl( X/S,
\log(1+p_n\o\CE^{\o n}) \bigr) \biggr) .
$$
\end{theorem}
\begin{proof}
By Theorem \ref{RESOLVE}, we know that
\begin{align*}
\Serre_S^{\SS_n}(\FF(X/S,n),j(n)^\bull\CE^{\boxtimes n})
&= \Serre_S^{\SS_n}(X^n/S,j(n)_\bull j(n)^\bull\CE^{\boxtimes n}) \\
&= \sum_{k=1}^n \Serre_S^{\SS_n}\biggl(X^n/S, \bigoplus_{J\in\S(n,k)}
\Hom(\s(n,J),\CE^{\boxtimes n}(J)) \biggr) \\
&= \sum_{k=1}^n \bigoplus_{J\in\S(n,k)}
\ch(\s(n,J)^\Dual)\o\Serre_S(X^n/S,\CE^{\boxtimes n}(J)) \in
\Lambda_n\o\K(S) .
\end{align*}
We may replace $\s(n,k)^\Dual$ by $\s(n,k)$, since any $\SS_n$-module is
isomorphic to its dual. Applying Theorem \ref{Lehrer-Solomon}, we obtain
$$
\sum_{k=1}^n \bigoplus_{J\in\S(n,k)} \prod_{i=1}^k \ch(\l(J_i)) \o
\Serre_S\bigl( X,\CE^{\o|J_i|} \bigr) ,
$$
where the product is taken in the ring $\Lambda\ohat\K(S)$. Summing over
$n$ gives
$$
\sum_{n=0}^\infty \Serre_S^{\SS_n}(X^n/S,j(n)_\bull
j(n)^\bull\CE^{\boxtimes n}) = \Exp \biggl( \sum_{n=1}^\infty
\ch(\l(n))\o\Serre_S(X,\CE^{\o n}) \biggr) .
$$
The theorem now follows from the character formula of Lemma
\ref{Hanlon-Stanley}, which may be rewritten as
$$
\ch(\l(n)) = \frac{1}{n} \sum_{d|n} (-1)^{n/d-1} \mu(d) p_d^{n/d} .
$$
We see that
\begin{align*}
\sum_{n=1}^\infty \l(n)\o\Serre_S(X,\CE^{\o n})) &= \sum_{n=1}^\infty
\frac{1}{n} \sum_{d|n} (-1)^{n/d-1} \mu(d) p_d^{n/d} \Serre_S(X,\CE^{\o n}) \\
&= \sum_{d=1}^\infty \frac{\mu(d)}{d} \sum_{e=1}^\infty
\frac{(-1)^{e-1}}{e} p_d^e \Serre_S(X,\CE^{\o de}) \\
&= \sum_{d=1}^\infty \frac{\mu(d)}{d} \Serre_S\bigl(X,\log(1+p_d\o\CE^{\o d})
\bigr) .
\qed\end{align*}
\def\qed{}
\end{proof}

\begin{remark}
Rewriting $\Exp$ in terms of Adams operations, Theorem \ref{push} becomes
$$
\sum_{n=0}^\infty \Serre_S^{\SS_n}(\FF(X/S,n),j(n)^\bull\CE^{\boxtimes n})
= \exp \biggl( \sum_{n=1}^\infty \frac{1}{n} \sum_{\ell=1}^\infty
\frac{(-1)^{\ell-1}}{\ell} p_n^\ell \sum_{d|n} \mu(n/d)
\psi_d\bigl(\Serre_S\bigl( X, \CE^{\o\ell n/d} \bigr) \bigr) \biggr) .
$$
For example, with the notation $\Serre(n)=\Serre_S(X,\CE^{\o n})$, we have
$$
\Serre_S^{\SS_n}(\FF(X/S,n),j(n)^\bull\CE^{\boxtimes n}) =
\sum_{\lambda\vdash n} s_\lambda\o \Phi_\lambda ,
$$
where $\Phi_{1^n}=\sigma_{1^n}(\Serre(1))$, while the other $\Phi_\lambda$
are as follows for $|\lambda|\le4$:
\begin{align*}
\Phi_2 &= \sigma_2(\Serre(1)) - \Serre(2) , \\
\Phi_3 &= \sigma_3(\Serre(1)) - \Serre(1)\Serre(2) , \quad
\Phi_{21} = \sigma_{21}(\Serre(1)) - \Serre(1)\Serre(2) + \Serre(3) , \\
\Phi_4 &= \sigma_4(\Serre(1)) - \sigma_2(\Serre(1))\Serre(2) +
\sigma_{1^2}(\Serre(2)) , \\
\Phi_{31} &= \sigma_{31}(\Serre(1)) -
\sigma_2(\Serre(1))\Serre(2) - \sigma_{11}(\Serre(1))\Serre(2) +
\Serre(1)\Serre(3) + \sigma_2(\Serre(2)) -\Serre(4) , \\
\Phi_{2^2} &= \sigma_{2^2}(\Serre(1)) - \sigma_2(\Serre(1))\Serre(2) +
\Serre(1)\Serre(3) + \sigma_{1^2}(\Serre(2)) , \\
\Phi_{21^2} &= \sigma_{21^2}(\Serre(1)) - \sigma_{1^2}(\Serre(1))\Serre(2)
+ \Serre(1)\Serre(3) - \Serre(4) .
\end{align*}
Note that the operations $\Phi_\lambda$ of \cite{I} are the specializations
of these polynomials obtained on setting $\Serre(n)=\Serre(1)$ for all
$n\ge1$.
\end{remark}

\medskip

Applying Theorem \ref{push} with $\CE$ equal to the unit object $\1$ of
$\D(X)$, and using that $\1^{\o n}=\1$ for all $n$, we obtain the following
corollary. Here, we abbreviate $\Serre_S(X,\1)$ to $\Serre_S(X)$.
\begin{corollary}
$\displaystyle \sum_{n=0}^\infty \Serre_S^{\SS_n}(\FF(X/S,n)) = \Exp \bigl(
\Log(1+p_1) \Serre_S(X) \bigr) $
\end{corollary}

Theorem \ref{push}, and its corollary, generalize immediately to the
equivariant situation, in which a finite group $\Gamma$ acts on $X$ and
$S$, and the morphism $\pi:X\to S$ and $\CE$ are $G$-equivariant. The
calculations now take place in the complete $\lambda$-ring
$\Lambda\ohat\K^\Gamma(S)$, and the formulas do not change.

\subsection{The configuration spaces of group schemes}
If $\G$ is a group scheme over $S$ and $n>0$, the scheme $\FF(\G/S,n)$ is
an $\SS_n$-equivariant $\G$-torsor, and we may consider the quotient scheme
$\G\backslash\FF(\G/S,n)$. Imitating the above proof, we now calculate its
$\SS_n$-equivariant relative Euler characteristic. There is also a
$\Gamma$-equivariant generalization, when a finite group $\Gamma$ acts on
all of the data; however, it is formally identical, so we simplify notation
by only treating the case $\Gamma=1$.
\begin{theorem} \label{relative-torsor}
If in the setting of Theorem \ref{push} $X=\G$ is a group scheme, then
$$
\sum_{n=1}^\infty \Serre_S^{\SS_n}(\G\backslash\FF(\G/S,n)) = \frac{\Exp
\bigl( \Log(1+p_1) \Serre_S(\G) \bigr) - 1}{\Serre_S(\G)} .
$$
\end{theorem}
\begin{proof}
We must first choose a $\G$-equivariant cover $\CU$ of $\FF(\G/S,2)$ by
affine open \embeddings. Observe that the automorphism
$(g,h)\mapsto(g,g^{-1}h)$ of $\G^2/S$ identifies $\FF(\G/S,2)$ with
$\G\times_S\G_0$, where $\G_0$ is the complement of the identity section of
$\G$. Under this identification, the action of $\G$ on $\FF(\G/S,2)$
corresponds to left translation in the first factor of $\G\times_S\G_0$.
We now choose a cover $\{U_i\}$ of $\G_0$ by affine open \embeddings
$j_i:U_i\hookrightarrow\G_0$; the cover $\CU$ of $\FF(\G/S,2)$ is the
pullback of this cover by the projection from
$\FF(\G/S,2)\cong\G\times_S\G_0$ to $\G_0$. (Here, we use that affine open
\embeddings are preserved under base change.)

If $n>0$, the morphism $j(n): \FF(\G/S,n) \hookrightarrow \G^n/S$ is an
$\SS_n$-equivariant \embedding of $\G$-torsors. On quotienting by the
action of the group scheme $\G$, we obtain an $\SS_n$-equivariant
\embedding
$$
\jbar(n): \G\backslash\FF(\G/S,n) \hookrightarrow \G\backslash(\G^n/S) .
$$
(Note that $\G\backslash \G^n/S$ is isomorphic to $\G^{n-1}/S$; however,
this isomorphism obscures the action of the symmetric group $\SS_n$.)
Since the cover $\CU$ of $\FF(\G/S,2)$ is $\G$-equivariant, the resolution
$\Resolve_\CU^\bull(\G/S,\1,n)$ is $\SS_n\times\G$-equivariant, so descends
to an $\SS_n$-equivariant resolution
$\G\backslash\Resolve_\CU^\bull(\G/S,\1,n)$ of
$\jbar(n)_\bull\jbar(n)^\bull\1$. The theorem now follows by a proof which
is entirely parallel to that of Theorem \ref{push} (in the special case
that $\CE=\1$), provided we observe that
\begin{align*}
\Serre_S^{\SS_n}(\G\backslash(\G^n/S),\G\backslash
\Resolve_\CU^{n-k}(\G/S,\1,n)) &= \sum_{k=1}^n \ch(\s(n,k))
\Serre_S(\G)^{k-1} \\ &= \frac{1}{\Serre_S(\G)}
\Serre_S^{\SS_n}(\G^n/S,\Resolve^{n-k}(\G/S,\1,n)) .
\qed\end{align*}
\def\qed{}
\end{proof}

Note that if $\G$ is a family of elliptic curves, the proof of Theorem
\ref{relative-torsor} simplifies, since we may take for the cover $\CU$ of
$\FF(\G/S,2)$ the canonical choice $\{\FF(\G/S,2)\}$. In fact, this is the
case of Theorem \ref{relative-torsor} which we apply in the next section,
to the universal family $E(N)$ of elliptic curves over the modular curve
$Y(N)$.

\section{The $\SS_n$-equivariant Serre \polynomial of the moduli space
$\CM_{1,n}$}

Let $\CM_{1,n}(N)$ be the fine moduli space of smooth elliptic curves of
level $N\ge3$ with $n$ marked points; it is a smooth quasi-projective
variety. The finite group $\SL(2,\Z/N)$ acts on $\CM_{1,n}(N)$, with
quotient $\CM_{1,n}$ the coarse moduli space of smooth elliptic curves.

Let $Y(N)$ be the modular curve $\CM_{1,1}(N)$, and let $E(N)\to Y(N)$ be
the universal elliptic curve of level $N$. The relative configuration space
$\FF(E(N)/Y(N),n)$ is an $\SL(2,\Z/N)$-equivariant $E(N)$-torsor with base
$Y(N)$.

Denote by $\HH$ the mixed Hodge module $\RR^1f_!\Q$ on $Y(N)$; it is of
course an $\SL(2,\Z/N)$-equivariant local system of rank $2$, known as the
Hodge local system. The sub-$\lambda$-ring which it generates in
$\K^{\SL(2,\Z/N)}(\MHM(Y(N)))$ is isomorphic to the Grothendieck group of
polynomial representations of the algebraic group $\GL(2)$; this is the
polynomial ring $\Z[\HH,\LL]$, with $\sigma_t(\HH)=(1-t\HH+t^2\LL)^{-1}$
and $\sigma_t(\LL)=(1-t\LL)^{-1}$. In this notation, we have
$$
\Serre_{Y(N)}^{\SL(2,\Z/N)}(E(N)) = 1 - \HH + \LL .
$$
We may now apply the $\SL(2,\Z/N)$-equivariant version of Theorem
\ref{relative-torsor}, obtaining the following formula.
\begin{proposition}
$$
\sum_{n=1}^\infty \Serre_{Y(N)}^{\SL(2,\Z/N)\times\SS_n}(\CM_{1,n}(N)) =
\frac{\displaystyle \biggl\{ \prod_{k=1}^\infty
(1+p_k)^{\frac{1}{k}\sum_{d|k}\mu(k/d) (1-\psi_d(\HH)+\LL^d)} \biggr\} - 1
} {1-\HH+\LL}
$$
\end{proposition}

Denote the $n$th symmetric power of $\HH$ by $\HH_n$; it is a rank $(n+1)$
$\SL(2,\Z/N)$-equivariant local system on $Y(N)$, given by the Chebyshev
polynomial of the second kind%
\footnote{These polynomials have generating function $\sum_{n=0}^\infty t^n
U_n(x) = (1-2xt+t^2)^{-1}$.}
$$
\HH_n = U_n(\HH/2) .
$$
The following table gives
$\Serre_{Y(N)}^{\SL(2,\Z/N)\times\SS_n}(\CM_{1,n}(N))$ for $n\le5$. This
table was calculated using J. Stembridge's symmetric function package
$\texttt{SF}$ \cite{SF} for $\mathtt{maple}$.
$$\begin{tabular}{|C|L|} \hline
n & \Serre_{Y(N)}^{\SL(2,\Z/N)\times\SS_n}(\CM_{1,n}(N)) \\ \hline
1 & \HH_0 \\[5pt]
2 & \HH_0 \LL s_2 - \HH_1 s_{1^2} \\[5pt]
3 & \HH_0 \LL^2 s_3 - \HH_1 (\LL s_{21} - s_{3}) + \HH_2 s_{1^3} \\[5pt]
4 & \HH_0 (\LL^3-\LL)s_4 - \HH_1(\LL^2s_{31}-\LL(s_4+s_{31})+s_{2^2}) +
\HH_2 (\LL s_{21^2} - s_{31}) - \HH_3 s_{1^4} \\[5pt]
5 & \HH_0 (\LL^4 s_5-\LL^2(s_5+s_{41})+\LL s_{32}) \\
& {} - \HH_1 (\LL^3
s_{41}-\LL^2(s_5+s_{41}+s_{32})-\LL(s_{32}+s_{2^21})+s_{31^2}) \\
& {}+ \HH_2 (\LL^2
s_{31^2}-\LL(s_{41}-s_{32}-s_{31^2})+(s_5+s_{32}+s_{2^21})) - \HH_3 (\LL
s_{21^3} - s_{31^2}) + \HH_4 s_{1^5} \\[3pt] \hline
\end{tabular}$$

\subsection{The Eichler-Shimura isomorphism}
Let $S_\ell(N)$ be the spaces of cusp forms of weight $\ell$ for the
congruence group $\Gamma(N) = \ker\bigl( \SL(2,\Z) \to \SL(2,\Z/N) \bigr)$.
It is an $\SL(2,\Z/N)$-module, and its invariant subspace
$S_\ell=S_\ell(1)$ is the space of cusp forms of level $1$.

Let $E_\ell(N)$ be the space of Eisenstein series of weight $\ell$. If
$\ell>2$, this is isomorphic as a $\SL(2,\Z/N)$-module to the induced
representation
$$
\Sigma_\ell(N) = \Ind^{\SL(2,\Z/N)}_{P(N)} \chi_\ell ,
$$
where $P(N)\subset\SL(2,\Z/N)$ is the parabolic subgroup of upper
triangular matrices, with generators
$T=\bigl[\begin{smallmatrix}1&1\\0&1\end{smallmatrix}\bigr]$ and $-I$, and
$\chi_\ell$ is the character of $P(N)$ which equals $1$ on $T$ and
$(-1)^\ell$ on $-I$.

The space $E_2(N)$ is smaller than $\Sigma_2(N)$: it is isomorphic to
$H_0(\SL(2,\Z/N),\Sigma_2(N)\bigr)$.

If $\ell$ is even, $\Sigma_\ell(N)$ is the permutation representation of
$\SL(2,\Z/N)$ on the set of cusps, and the $\SL(2,\Z/N)$-invariant subspace
is one-dimensional; the corresponding subspace of $E_\ell(N)$ is spanned by
the level $1$ Eisenstein series $E_\ell$. If $\ell$ is odd, there are no
$\SL(2,\Z/N)$-invariant elements of $E_\ell(N)$, reflecting the fact that
there are no level $1$ Eisenstein series of odd weight. In all cases,
$\Sigma_\ell(N)$ has dimension
$$
[\SL(2,\Z/N):P(N)] = \frac{k^2}{2} \prod_{p|k}(1-p^{-2})
$$
equal to the number of cusps of the congruence subgroup $\Gamma(N)$.

Eichler and Shimura have calculated the cohomology of the sheaves $\HH_n$.
This calculation is explained in Verdier \cite{Verdier} and Shimura
\cite{Shimura}. The mixed Hodge structure on this cohomology may be
calculated by the same technique that Deligne uses in \cite{Deligne} to
calculate the action of the Frobenius operator on the \'etale cohomology
groups. Define the Hodge structure $\SSS_{n+2}(N)$ to be
$\gr^W_{n+1}H^1_c(Y(N),\HH_n)$.
\begin{theorem} \label{Shimura}
The vector spaces $\gr_k^WH^i_c(Y(N),\HH_n)$ associated to the weight
filtration on the cohomology groups $H^\bull_c(Y(N),\HH_n)$ vanish, with
the exception of
\begin{align*}
& \gr^W_0H^1_c(Y(N),\HH_n)\cong E_{n+2}(N) , \\
& \gr^W_{n+1}H^1_c(Y(N),\HH_n)\cong \SSS_{n+2}(N) , \quad\text{and} \\
& \gr^W_2H^2_c(Y(N),\HH_0)=\LL .
\end{align*}
The Hodge filtration of $\SSS_{n+2}(N)$ has two steps: $0\subset
F^0\SSS_{n+2}(N)\subset \SSS_{n+2}(N)$, and the vector space
$F^0\SSS_{n+2}(N)$ is naturally isomorphic to $S_{n+2}(N)$.
\end{theorem}

\begin{corollary} \label{final}
The equivariant Serre \polynomial
$$
\Serre^{\SL(2,\Z/N)\times\SS_n}(\CM_{1,n}(N)) \in
\K^{\SL(2,\Z/N)}(\MHM(\point)) \o \Lambda_n
$$
is obtained from the equivariant Serre \polynomial
$$
\Serre_{Y(N)}^{\SL(2,\Z/N)\times\SS_n}(\CM_{1,n}(N)) \in
\K^{\SL(2,\Z/N)}(\MHM(Y(N))) \o \Lambda_n
$$
by the substitution
$\HH_n\mapsto\delta_{n,0}(\LL+1)-\Sigma_{n+2}(N)-\SSS_{n+2}(N)$.
\end{corollary}

We may now descend to level $1$ by applying the augmentation
$$
\eps : \K^{\SL(2,\Z/N)}(\MHM(\point)) \to \K(\MHM(\point)) ,
$$
given explicitly by $\eps(\SSS_\ell(N))=\SSS_\ell$ and
$$
\eps(\Sigma_\ell(N)) = \begin{cases} \Q , & \text{$\ell$ even,} \\
  0 , & \text{$\ell$ odd.} \end{cases}
$$
The following table gives the $\SS_n$-equivariant Serre \polynomial
$\Serre^{\SS_n}(\CM_{1,n})$ for $n\le5$, together with the underlying Serre
\polynomial $\Serre(\CM_{1,n})$ and Euler characteristic.
$$\begin{tabular}{|C|L|L|L|} \hline
n & \Serre^{\SS_n}(\CM_{1,n}) & \Serre(\CM_{1,n}) & \chi(\CM_{1,n}) \\ \hline
1 & \LL s_1 & \LL & 1 \\[5pt]
2 & s_2 \LL^2 & \LL^2 & 1 \\[5pt]
3 & s_3 \LL^3 - s_{1^3} & \LL^3-1 & 0 \\[5pt]
4 & s_4 \LL^4 - s_4 \LL^2 - s_{21^2} \LL + s_{31} &
\LL^4 - \LL^2 - 3\,\LL + 3 & 0 \\[5pt]
5 & s_5 \LL^5 - (s_5+s_{41})\LL^3 + (s_{32}-s_{31^2})\LL^2
& \LL^5 - 5\,\LL^3 - \LL^2 + 15\,\LL - 12 & -2\\
& {}+ (s_{41}+s_{32}+s_{31^2})\LL & & \\
& {}- (s_5+s_{32}+s_{2^21}+s_{1^5}) & & \\[3pt] \hline
\end{tabular}$$

Corollary \ref{final} may be expressed in closed form:
\begin{multline} \label{closed}
\sum_{n=1}^\infty \Serre^{\SS_n}(\CM_{1,n}) = \Res_0 \Biggl[ \left(
\frac{\prod_{n=1}^\infty
(1+p_n)^{\frac{1}{n}\sum_{d|n}\mu(n/d)(1-\om^d-\LL^d/\om^d+\LL^d)} - 1}
{1-\om-\LL/\om+\LL} \right) \\ \times \left( \sum_{k=1}^\infty \biggl(
\frac{\SSS_{2k+2}+1}{\LL^{2k+1}} \biggr) \om^{2k} - 1 \right) \bigl(
\om-\LL/\om \bigr) d\om \Biggr] ,
\end{multline}
where $\Res_0$ is the residue of the differential form at the origin. This
is an easy consequence of the Weyl integration formula for $\SU(2)$, in the
form
$$
- \frac{1}{2} \Res_0\biggl[ \biggl(
\frac{\om^{k+1}-(\LL/\om)^{k+1}}{\om-\LL/\om} \biggr) \biggl(
\frac{\om^{\ell+1}-(\LL/\om)^{\ell+1}}{\om-\LL/\om} \biggr) \bigl(
\om-\LL/\om \bigr)^2 \frac{d\om}{\om} \biggr] = \LL^{k+1} \delta_{k\ell} .
$$

To obtain a formula for the non-equivariant Serre \polynomials, we replace
$p_n$, $n>1$, by $0$, and expand in $p_1$, which gives
\begin{equation} \label{non-equi}
\frac{\Serre(\CM_{1,n+1})}{n!} = \Res_0 \Biggl[ \binom{\LL-\om-\LL/\om}{n}
\left( \sum_{k=1}^\infty \biggl( \frac{\SSS_{2k+2}+1}{\LL^{2k+1}} \biggr)
\om^{2k} - 1 \right) \bigl( \om-\LL/\om \bigr) d\om \Biggr] .
\end{equation}

{}From \eqref{non-equi}, we can calculate the Euler characteristic
$\chi(\CM_{1,n})$ directly. The following proof was shown to the author by
D. Zagier.
\begin{proposition} \label{don}
If $n\ge5$, $\chi(\CM_{1,n})=(-1)^n(n-1)!/12$.
\end{proposition}
\begin{proof}
If in \eqref{non-equi}, we replace $\LL=1$ and $\SSS_{k+2}$ by
$2\dim(S_{k+2})$, we see that
$$
\frac{\chi(\CM_{1,n+1})}{n!} = \Res_0 \biggl[ \binom{1-\om-\om^{-1}}{n}
\frac{(1-\om^2-2\om^4-\om^6+\om^8)}{(1+\om^2)(1-\om^6)} \frac{d\om}{\om}
\biggr] .
$$
The poles of this differential form are all simple, and are located at
$\om=0$ and $\om=\infty$, and at values of $\om$ such that $\om+\om^{-1}$
is an integer in the interval $[-2,2]$ (the latter poles are on the unit
circle). Since it is invariant under $\om\mapsto\om^{-1}$, its residues at
$0$ and $\infty$ are equal. By the residue theorem, it follows that
$$
\frac{\chi(\CM_{1,n+1})}{n!} = - \frac12 \sum_{z\in\{\pm1,\pm
i,\pm\rho,\pm\rho^2\}} \Res_z \biggl[ \binom{1-\om-\om^{-1}}{n}
\frac{(1-\om^2-2\om^4-\om^6+\om^8)}{(1+\om^2)(1-\om^6)} \frac{d\om}{\om}
\biggr] ,
$$
where $\rho$ is a primitive sixth root of unity.

The residues of this differential form on the unit circle are as follows:
\begin{equation}\label{residues}
\Res_z\biggl[ \dfrac{(1-\om^2-2\om^4-\om^6+\om^8)}{(1+\om^2)(1-\om^6)}
\biggr] = \begin{cases} 1/6 , & |z+z^{-1}| = 2 , \\
-1/3 , & |z+z^{-1}|=1 , \\
-1/2 , & |z+z^{-1}|=0 .
\end{cases}
\end{equation}
At each of these poles except $\om=1$, the binomial coefficient
$\binom{1-\om-\om^{-1}}{n}$ vanishes for $n\ge4$. This leaves the residue
at $1$, which equals $(-1)^{n+1}/12$.
\end{proof}

We close the paper with a calculation of the Serre \polynomials of the
spaces $\CM_{1,n}/\SS_n$. If we substitute $x^n$ for $p_n$ in
\eqref{closed}, we obtain the generating function for the $\SS_n$-invariant
parts of the local systems $\Serre^{\SS_n}(\CM_{1,n}/\CM_{1,1})$. By
Corollary (5.7) of \cite{I}, we have
\begin{multline*}
\sum_{n=1}^\infty H^0\bigl( \SS_n,\Serre^{\SS_n}(\CM_{1,n}/\CM_{1,1})
\bigr) x^n = \frac{\displaystyle\prod_{n=1}^\infty
(1+x^n)^{\frac{1}{n}\sum_{d|n}\mu(n/d)(1-\om^d-\LL^d/\om^d+\LL^d)} - 1}
{1-\om-\LL/\om+\LL} \\
\begin{aligned}
{} &= \frac{1}{1-\om-\LL/\om+\LL} \left\{ \frac{(1-\om x) (1-\LL x/\om)
(1-x^2)(1-\LL x^2)} {(1-x)(1-\LL x)(1-\om x^2)(1-\LL x^2/\om)} - 1 \right\}
\\
{} & = x \left( \frac{1-\LL x^3}{1-\LL x} \right)
\frac{1}{1-(\om+\LL/\om)x^2+\LL x^4} \\
{} & = x \left( \frac{1-\LL x^3}{1-\LL x} \right) \sum_{k=0}^\infty \HH_k
x^{2k} .
\end{aligned}
\end{multline*}
Applying the functor $H^\bull_c(\CM_{1,1},-)$, we see that:
$$
\sum_{n=1}^\infty \Serre(\CM_{1,n}/\SS_n) x^n = x \left( \frac{1-\LL x^3}
{1-\LL x} \right) \sum_{k=0}^\infty \Serre(\CM_{1,1},\HH_{2k}) x^{4k}
$$
Taking Euler characteristics gives
\begin{align*}
\sum_{n=1}^\infty \chi(\CM_{1,n}/\SS_n) x^n &=(x+x^2+x^3) \sum_{n=0}^\infty
\chi(\CM_{1,1},\HH_n) x^{4n} \\
& = (x+x^2+x^3) \frac{(1-x^4-2x^8-x^{12}+x^{16})}{(1-x^8)(1-x^{12})} .
\end{align*}
The corresponding formulas in genus $0$ are
$\Serre(\CM_{0,n}/\SS_n)=\LL^{n-3}$ and $\chi(\CM_{0,n}/{\SS_n})=1$, for
all $n\ge3$.

\end{document}